\documentclass[12pt]{article}                
\usepackage{graphicx}
\usepackage{psfrag}
\usepackage{amsmath, amssymb}
\usepackage{mathptmx}
\usepackage{ntheorem}

\makeatletter
\newsavebox\myboxA
\newsavebox\myboxB
\newlength\mylenA
\newcommand*\xoverline[2][0.75]{%
    \sbox{\myboxA}{$\m@th#2$}%
    \setbox\myboxB\null
    \ht\myboxB=\ht\myboxA%
    \dp\myboxB=\dp\myboxA%
    \wd\myboxB=#1\wd\myboxA
    \sbox\myboxB{$\m@th\overline{\copy\myboxB}$}
    \setlength\mylenA{\the\wd\myboxA}
    \addtolength\mylenA{-\the\wd\myboxB}%
    \ifdim\wd\myboxB<\wd\myboxA%
       \rlap{\hskip 0.5\mylenA\usebox\myboxB}{\usebox\myboxA}%
    \else
        \hskip -0.5\mylenA\rlap{\usebox\myboxA}{\hskip 0.5\mylenA\usebox\myboxB}%
    \fi}
\makeatother

\newcommand{\Av}{{\rm Av}}

\newcommand{\PP}{{\cal P}}

\newcommand{\eps}{{\varepsilon}}

\newcommand{\e}{\mathbb{E}}

\newcommand{\Reals}{\mathbb{R}}

\newcommand{\la}{\langle}
\newcommand{\ra}{\rangle}

\newtheorem{theorem}{Theorem}

\theoremstyle{nonumberplain}
\newcommand\specialref{}

\begin{document}

\title{On the $K$-sat model with large number of clauses
}
\author{Dmitry Panchenko\thanks{\textsc{\tiny Department of Mathematics, University of Toronto, panchenk@math.toronto.edu. Partially supported by NSERC.}
}\\
}
\date{}
\maketitle
\begin{abstract}
We show that in the $K$-sat model with $N$ variables and $\alpha N$ clauses, the expected ratio of the smallest number of unsatisfied clauses to the number of variables is $\alpha/2^K - \sqrt{\alpha} c_*(N)/2^K$ up to smaller order terms $o(\sqrt{\alpha})$ as $\alpha\to\infty$ uniformly in $N$, where $c_*(N)$ is the expected normalized maximum energy of some specific mixed $p$-spin spin glass model. The formula for the limit of $c_*(N)$ is well known in the theory of spin glasses.
\end{abstract} 
\vspace{0.5cm}

\section{Introduction}

Let $K\geq 2$ be an integer and let $\alpha>0$. Given $N\geq 1$, we will denote the elements of the hypercube $\{-1,+1\}^N$ by $\sigma=(\sigma_1,\ldots,\sigma_N)$. Consider an i.i.d. sequence of indices $(i_{j,k})_{j,k\geq 1}$ with uniform distribution on $\{1,\ldots,N\}$ and let $\pi(\alpha N)$ be an independent Poisson random variable with the mean $\alpha N.$ We define the $K$-sat Hamiltonian $H_\alpha(\sigma)$ on $\{-1,+1\}^N$ by
\begin{equation} 
H_{\alpha}(\sigma)=-\sum_{j\leq \pi(\alpha N)} \prod_{k\leq K} \frac{1+\eps_{j,k} \sigma_{i_{j,k}}}{2},
\label{Halpha}
\end{equation}
where $\eps_{j,k}$ are random signs (symmetric $\{-1,+1\}$-valued random variables) independent over different indices $(j,k)$ and independent of all other random variables. Each random clause 
$$
\prod_{k\leq K} \frac{1+\eps_{j,k} \sigma_{i_{j,k}}}{2}
$$
in the above sum can take values $0$ or $1$, depending on the values of the coordinates (Boolean variables) $ \sigma_{i_{j,k}}\in\{-1,+1\}$ participating in the clause. If the value is zero, the clause is said to be satisfied and, otherwise, it is unsatisfied. In other words, the clause represents a random disjunction because it is satisfied if at least one $\sigma_{i_{j,k}}=-\eps_{j,k}.$ With this convention, the quantity
\begin{equation}
M_{N,\alpha} = \max_{\sigma}\frac{H_{\alpha}(\sigma)}{N} = - \min_{\sigma}\frac{-H_{\alpha}(\sigma)}{N}
\end{equation}
represents (up to the minus sign) the smallest proportion of unsatisfied clauses over all possible assignments of $\sigma.$ For any $K\geq 3,$ it is expected that, up to a certain threshold, for $\alpha\leq \alpha_K$, with high probability all clauses can be satisfies and $\max_\sigma H_{\alpha}(\sigma) = 0,$ while above this threshold with high probability all clauses  can not be satisfied. The value of $\alpha_K$ was described precisely (in the sense of theoretical physics) by Mertens, M\'ezard and Zecchina in \cite{MMZ} on the basis of the celebrated M\'ezard-Parisi ansatz \cite{Mezard}, further developed  in \cite{MPZ}. For example, for $K=3$ the phase transition was predicted to be at $\alpha_3 \approx 4.267$, and the large $K$ behaviour to be 
\begin{equation}
\alpha_K=2^K \ln 2-\frac{1}{2}(1+\ln 2)+ o_K(1).
\label{thresh}
\end{equation}
This problem has been studied extensively in the mathematics literature, with progressively more precise results obtained in \cite{AM02, AP03, CP13, COKSAT}, and the exact threshold for large enough $K$ was finally determined in \cite{DSS}. Describing the threshold for all $K\geq 3$ remains an open problem.

In this paper we will consider the regime of large $\alpha$, in which case the proportion of unsatisfied clauses is strictly positive. Given $\sigma$, if we select a clause randomly, the probability of it being unsatisfied is $1/2^K$. It turns out that, for large $\alpha$, optimal assignments are not much better than any fixed assignment and the leading term of the smallest ratio of unsatisfied clauses to the number of variables $N$ is $\alpha/2^K.$ We will show that for optimal assignments the next order correction term for large $\alpha$ is of the form $-c_* \sqrt{\alpha}$, where the constant $c_*=c_*(N)$ is related to the expected maximum of the specific mixed $p$-spin spin glass model in (\ref{H}) below. This will establish the Leuzzi-Parisi formula obtained in \cite{LP} by the non-rigorous replica method. Let us consider the following mixed $p$-spin Hamiltonian 
\begin{equation}
H(\sigma) = \sum_{p=1}^K  \sqrt{{K\choose p} \frac{1}{N^{p-1}}}
\sum_{1\leq i_1,\ldots,i_p\leq N}g_{i_1,\ldots,i_p} \sigma_{i_1}\cdots\sigma_{i_p},
\label{H}
\end{equation}
where the coefficients $(g_{i_1,\ldots,i_p})$ are standard Gaussian random variables independent for all $p\geq 1$ and all indices $(i_1,\ldots,i_p)$. If we consider the function
\begin{equation}
\xi(x) = \sum_{p=1}^K {K\choose p} x^p = (1+x)^K-1
\label{xi}
\end{equation}
then the covariance of the Gaussian Hamiltonian (\ref{H}) is given by
\begin{equation}
\e H(\sigma^1) H(\sigma^2) = N \xi(R_{1,2}), 
\end{equation}
where 
\begin{equation}
R_{1,2}= \frac{1}{N}\sum_{i\leq N} \sigma^1_i \sigma_i^{2}
\label{overlap}
\end{equation}
is the overlap of configurations $\sigma^1$ and  $\sigma^2$. Let us denote the normalized maximum by
\begin{equation}
M_N = \max_{\sigma}\frac{H(\sigma)}{N}.
\end{equation}
Our main result is the following.
\begin{theorem}\label{Th1}
For all $N\geq 1$, we have
\begin{equation}
\e M_{N,\alpha} = - \frac{\alpha}{2^K} + \frac{\sqrt{\alpha}}{2^K} \e M_N + R(\alpha),
\end{equation}
where $|R(\alpha)|\leq L\alpha^{1/3}$ for $\alpha\geq L$ for some absolute constant $L.$
\end{theorem}
Notice that the remainder term is guaranteed to be smaller than the correction term only when $\alpha$ is of the order $(2^K)^6$, which is way above the phase transition (\ref{thresh}) for large $K.$ The proof of Theorem \ref{Th1} is based on the interpolation technique of Guerra and Toninelli in \cite{GT} (not to be confused with another Guerra-Toninelli interpolation \cite{GuerraToninelli}) and is a slight modification of the argument in \cite{DMS}, were similar results for extremal cuts of sparse random graphs were obtained. More recent results in this direction, for example, for diluted $p$-spin spin glass models, can be found in \cite{Sen}. Further applications of the Guerra-Toninelli interpolation \cite{GT} can be found in \cite{JKS, ChenP17}.

Perhaps, the main reason why Theorem \ref{Th1} is interesting is because mixed $p$-spins models are much better understood than diluted models and, in particular, the formula for the limit of $\e M_N$ is known, while previously only upper and lower bounds on the factor in front of $\sqrt{\alpha}$ were known (see Theorem 15 in \cite{CGHS}). This limit can be expressed as the zero temperature limit of the celebrated Parisi formula \cite{Parisi79, Parisi} for the free energy of the mixed $p$-spin models. The first proof of the Parisi formula for mixed even $p$-spin models was obtained by Talagrand in \cite{TPF}, building upon the replica symmetry breaking interpolation method of Guerra \cite{Guerra}. The model we consider in (\ref{H}) includes odd $p$-spin interaction terms and in this generality the Parisi formula was proved in \cite{PPF} as a consequence of the Parisi ultrametricity hypothesis for the overlaps proved in \cite{PUltra} (see also \cite{SKmodel}). 

The good news is that, due to a recent breakthrough in \cite{ChenAuf} (building upon the ideas in \cite{ChenAufP}), the zero temperature limit of the Parisi formula can be expressed in a form (conjectured by Guerra) quite similar to the classical Parisi formula at positive temperature, as follows. Let $\cal U$ be the family of all nonnegative nondecreasing step functions on $[0,1]$ with finitely many jumps. For $u\in {\cal U}$, let $\Psi_u(t,x)$ for $(t,x)\in[0,1]\times\Reals$ be the solution of
\begin{equation}
\frac{\partial \Psi_u}{\partial t} = -\frac{1}{2}\xi''(t)
\Bigl(
\frac{\partial^2 \Psi_u}{\partial x^2} + u(t)
\Bigl(\frac{\partial \Psi_u}{\partial x}\Bigr)^2
\Bigr)
\label{heat}
\end{equation}
with the boundary condition $\Psi_u(1,x)=|x|.$ Define 
\begin{equation}
\PP(u) = \Psi_u(0,0) - \frac{1}{2}\int_0^1 t\xi''(t)u(t)\, dt.
\end{equation}
Then Theorem 1 in \cite{ChenAuf} shows that 
\begin{equation}
\lim_{N\to\infty} \e M_N = \inf_{u\in {\cal U}} \PP(u).
\end{equation}
We refer to \cite{ChenAuf} for further details and turn to the proof of Theorem \ref{Th1}.

\section{Proof of the main result}

For $t\in [0,1]$, let us consider the interpolating Hamiltonian
\begin{equation}
H(t,\sigma) = \delta H_{\alpha(1-t)}(\sigma) +\sqrt{t}\beta H(\sigma),
\label{Hampert}
\end{equation}
where the first term $H_{\alpha(1-t)}(\sigma)$ is defined as in (\ref{Halpha}), only with $\alpha$ replaced by $\alpha(1-t),$ and the inverse temperature parameters $\delta>0$ and $\beta>0$ will be chosen later. Let
\begin{equation}
\varphi(t) = \frac{1}{N}\e\log \sum_{\sigma} \exp H(t,\sigma)
\end{equation}
be the corresponding interpolating free energy. It is a well-known and straightforward calculation to compute the derivative $\varphi'(t)$ using Gaussian integration by parts for the second term and Poisson integration by parts for the first term in (\ref{Hampert}). It can be written as $\varphi'(t)=\mathrm{I}+\mathrm{II}$, with the two terms defined as follows. Let us denote by $\la\,\cdot\,\ra_t$ the average with respect to the Gibbs measure
$$
G_t(\sigma) = \frac{\exp H(t,\sigma)}{\sum_\sigma \exp H(t,\sigma)}
$$
corresponding to the Hamiltonian $H(t,\sigma)$, as well as the average with respect to its infinite product  $G_t^{\otimes\infty}$. Taking the derivative in $\sqrt{t}$ in the second term in (\ref{Hampert}) and using standard Gaussian integration by parts (see e.g. \cite{SG2} or Section 1.2 in \cite{SKmodel}),
$$
\mathrm{I} = \frac{\beta^2}{2}\Bigl(\xi(1)- \e\bigl\la\xi(R_{1,2})\bigr\ra_t\Bigr),
$$
where $\xi$ was defined in (\ref{xi}) and $R_{1,2}$ is the overlap in (\ref{overlap}). To write the second term $\mathrm{II}$, let us introduce the notation, for $n\geq 1$,
\begin{equation}
Q_{1,\ldots,n} =  \frac{1}{N}\sum_{i\leq N} \Av \prod_{\ell \leq n}\frac{1+\eps\sigma_i^\ell}{2},
\label{multioverlaps}
\end{equation}
where $\Av$ is the average over $\eps=\pm 1$ with equal weights $1/2.$ For example,
\begin{equation}
Q_{1} = \frac{1}{2} \,\mbox{ and }\, Q_{1,2} = \frac{1+R_{1,2}}{4}.
\end{equation}
Then, a standard argument using Poisson integration by parts (see e.g. \cite{FL}, \cite{PT} or \cite{Pspins}) gives,
$$
\mathrm{II} = \alpha
\sum_{n\geq 1}\frac{(1-e^{-\delta})^n}{n}
\e \bigl\la (Q_{1,\ldots,n})^K\bigr\ra_t.
$$
The first two terms on the right hand side are equal to
\begin{align*}
& \frac{\alpha(1-e^{-\delta})}{2^K}+ \frac{\alpha(1-e^{-\delta})^2}{2\cdot 4^K} \e \bigl\la (1+R_{1,2})^K\bigr\ra_t
\\
& =
\frac{\alpha(1-e^{-\delta})}{2^K}+ \frac{\alpha(1-e^{-\delta})^2}{2\cdot 4^K} 
+ \frac{\alpha(1-e^{-\delta})^2}{2\cdot 4^K}  \e \bigl\la \xi(R_{1,2})\bigr\ra_t,
\end{align*}
and we will denote the remainder by
\begin{equation}
\mathrm{III} = \alpha \sum_{n\geq 3}\frac{(1-e^{-\delta})^n}{n} \e \bigl\la (Q_{1,\ldots,n})^K\bigr\ra_t.
\label{Remainder}
\end{equation}
For a given $\alpha$ and $\delta$, we are going to make the following choice of $\beta$,
\begin{equation}
\frac{\beta^2}{2} = \frac{{\alpha}(1-e^{-\delta})^2}{2\cdot 4^K}, \,\mbox{ or }\,
\beta = \frac{\sqrt{\alpha}(1-e^{-\delta})}{2^K}.
\end{equation}
With this choice, the coefficients in front of $\e \la \xi(R_{1,2})\ra_t$ in the terms I and II cancel out and, using that $\xi(1)+1=2^K$, we can express
\begin{align}
\varphi'(t) 
&= \frac{\beta^2}{2}\xi(1) + \frac{\alpha(1-e^{-\delta})}{2^K}+ \frac{\alpha(1-e^{-\delta})^2}{2\cdot 4^K} +\mathrm{III} 
\nonumber
\\
& = \frac{\alpha(1-e^{-\delta})}{2^K}+ \frac{\alpha(1-e^{-\delta})^2}{2\cdot 2^K} +\mathrm{III}.
\label{phider}
\end{align}
The rest of the proof is a collection of elementary estimates. First of all, if we denote $x = 1-e^{-\delta}$ then $-\log(1-x)=\delta$ and, since $Q_{1,\ldots,n}\in [0,1]$,
$$
|\mathrm{III}|\leq \alpha \sum_{n\geq 3}\frac{(1-e^{-\delta})^n}{n}
= \alpha\Bigl( -\log(1-x) - x - \frac{x^2}{2}\Bigr) = O(\alpha x^3) = O(\alpha \delta^3)
$$
for $\delta$ small enough. Integrating (\ref{phider}) between $0$ and $1$, we get
$$
\Bigl| \varphi(0) + \frac{\alpha(1-e^{-\delta})}{2^K}+ \frac{\alpha(1-e^{-\delta})^2}{2\cdot 2^K} - \varphi(1)\Bigr| = O(\alpha \delta^3)
$$
and, dividing both sides by $\delta$,
$$
\Bigl| \frac{1}{\delta} \varphi(0) + \frac{\alpha(1-e^{-\delta})}{2^K\delta}+ \frac{\alpha(1-e^{-\delta})^2}{2\cdot 2^K\delta} - \frac{1}{\delta}\varphi(1)\Bigr| = O(\alpha \delta^2).
$$
Next, we will need the following estimates,
$$
\frac{1}{N}\e\max_\sigma H(t,\sigma) \leq \varphi(t) = \frac{1}{N}\e\log \sum_{\sigma} \exp H(t,\sigma)
\leq \log 2 + \frac{1}{N}\e\max_\sigma H(t,\sigma),
$$
which can be obtained by keeping only the largest term in the sum $\sum_\sigma$ in the middle to get the lower bound, and replacing all $2^N$ terms by the largest one to get the upper bound. Using this for $t=0$ and $t=1$, we get
$$
\Bigl| \frac{1}{\delta} \varphi(0) - \e M_{N,\alpha}\Bigr| \leq \frac{\log 2}{\delta},\,\,
\Bigl| \frac{1}{\delta} \varphi(1) -  \frac{\beta}{\delta}\e M_N \Bigr| \leq \frac{\log 2}{\delta}.
$$
Therefore,
$$
\Bigl| \e M_{N,\alpha} + \frac{\alpha(1-e^{-\delta})}{2^K\delta}+ \frac{\alpha(1-e^{-\delta})^2}{2\cdot 2^K\delta} - \frac{\sqrt{\alpha}(1-e^{-\delta})}{2^K\delta} \e M_N \Bigr| = O\Bigl(\frac{1}{\delta}+ \alpha \delta^2\Bigr).
$$
By Taylor's expansion, for $\delta$ small enough,
$$
\frac{\alpha(1-e^{-\delta})}{2^K\delta}+ \frac{\alpha(1-e^{-\delta})^2}{2\cdot 2^K\delta} = \frac{\alpha}{2^K} + O(\alpha \delta^2)
$$
and
$$
\frac{\sqrt{\alpha}(1-e^{-\delta})}{2^K\delta} = \frac{\sqrt{\alpha}}{2^K} +O(\sqrt{\alpha\delta^2}).
$$
Plugging this in the above equation,
$$
\Bigl| \e M_{N,\alpha} + \frac{\alpha}{2^K} - \frac{\sqrt{\alpha}}{2^K} \e M_N\Bigr| = O\Bigl(\frac{1}{\delta}+ \alpha \delta^2+\sqrt{\alpha\delta^2}\Bigr).
$$
Taking $\delta=\alpha^{-1/3}$ finishes the proof.

\section{Some comments}

We saw in the above proof that, in order to obtain the first correction term, we could discard and roughly bound all the remainder terms in (\ref{Remainder}) and match only the second term involving the overlap $R_{1,2}$ with the corresponding Gaussian mixed $p$-spin model. In order to compute the limit $\lim_{N\to\infty} \e M_{N,\alpha}$ precisely for any fixed $\alpha$, one approach is to take the zero temperature limit of the M\'ezard-Parisi formula \cite{Mezard} for the free energy of diluted models at positive temperature and it seems that, in order to prove this formula rigorously, understanding the remainder terms is crucial. The M\'ezard-Parisi formula can be derived using some known results if one can figure out how to prove the key hypothesis made in \cite{Mezard}, namely, that multi-overlaps $Q_{1,\ldots,n}$ in (\ref{multioverlaps}) are continuous functions of the overlaps $R_{\ell,\ell'}$ for $\ell,\ell'\leq n$. So far, this was shown in \cite{finiteRSB} (in some perturbative sense) only under the technical assumption that, asymptotically, the Gibbs measure of the model has finite many steps of replica symmetry breaking or, in other words, the overlap $R_{1,2}$ takes only finitely many values. Proving this hypothesis in full generality remains one of the main obstacles in understanding the $K$-sat and other diluted models for any fixed $\alpha.$

\end{document}